\documentclass[12pt]{amsart}
\usepackage{amscd,amssymb}

\let\oldlabel=\label
\def\prellabel{\marginparsep=1em\marginparwidth=44pt
    \def\label##1{\oldlabel{##1}\ifmmode\else\ifinner\else
         \marginpar{{\footnotesize\ \\ \tt
                    ##1}}\fi\fi}}
\def\NN{{\mathbb N}}
\def\op{\operatorname{op}}
\def\W{{\operatorname{Witt}}}
\def\vert{\operatorname{vert}}
\def\X{\operatorname{X}}
\def\gp{\operatorname{gp}}
\def\rank{\operatorname{rank}}

\def\Coker{\operatorname{Coker}}
\def\Proj{\operatorname{Proj}}
\def\chara{\operatorname{char}}

\def\Proj{\operatorname{Proj}}
\def\Spec{\operatorname{Spec}}
\def\RR{{\mathbb R}}
\def\QQ{{\mathbb Q}}
\def\ZZ{{\mathbb Z}}
\def\NN{{\mathbb N}}
\def\TT{{\mathbb T}}
\def\AA{{\mathbb A}}
\def\Zz{{\mathcal Z}}
\def\Aa{{\mathcal A}}
\def\Ff{{\mathcal F}}
\def\C{{\mathfrak C}}

\advance\headheight1.15pt

\newtheorem{lemma}{Lemma}%[section]
\newtheorem{corollary}[lemma]{Corollary}
\newtheorem{theorem}[lemma]{Theorem}
\newtheorem{proposition}[lemma]{Proposition}

\theoremstyle{definition}

\newtheorem{example}[lemma]{Example}

\begin{document}

\title[Toric varieties with huge Grothendieck group]
{Toric varieties with huge Grothendieck group}

\author{Joseph Gubeladze}

\address{Department of Mathematics, San Francisco state University, San
Francisco, CA 94132, USA}

\address{A. Razmadze Mathematical Institute, Alexidze St. 1, 380093
Tbilisi, Georgia}

\email{soso@math.sfsu.edu}

\thanks{The work was done during the author's visit to Universit\'e
Louis Pasteur, Strasbourg (Autumn 2001). He was also supported by
INTAS grant 99-00817 and TMR grant ERB FMRX CT-97-0107}


\begin{abstract}
In each dimension $n\geq3$ there are many projective simplicial
toric varieties whose Grothendieck groups of vector bundles are at
least as big as the ground field. In particular, the conjecture
that the Grothendieck groups of locally trivial sheaves and
coherent sheaves on such varieties are rationally isomorphic fails
badly.
\end{abstract}

\maketitle

\section{Introduction}

The question whether the natural homomorphism $K_0(X)\to G_0(X)$
between the Grothendieck groups of locally trivial sheaves and
coherent sheaves on a complete simplicial toric variety $X$ is a
rational isomorphism  has attracted an attention of researchers
for some time (see \cite{BV}). More generally, it has been
conjectured that such a map is an isomorphism after tensoring with
$\QQ$ for arbitrary quasi-projective orbifold (see, for instance,
\cite[Section 7.12]{C}). The main result of this paper says that
this is far from being the case even for \emph{projective
simplicial toric varieties} in dimensions $3$ and higher.

More precisely, it was shown in \cite{BV} that $K_0(X)_{\QQ}\to
G_0(X)_{\QQ}$ is a surjection for a simplicial toric variety $X$.
If $X$ is quasiprojective then by {\it Riemann-Roch for singular
varieties} \cite{BFM} we have an isomorphism $G_0(X)_{\QQ}\to
A_*(X)_{\QQ}$. On the other hand the Chow groups of $X$ are
finitely generated by \cite{FMSS}. Therefore, $G_0(X)$ has a
finite rank. Fix arbitrary field $k$, not necessarily
algebraically closed and satisfying the conditions $\chara k=0$
and $\dim_{\QQ}k=\infty$. Theorem \ref{main} below gives many
examples of projective simplicial toric $k$-varieties $X$ such
that $\rank K_0(X)\geq\dim_{\QQ}k$.

A word on notation and terminology. Throughout the paper $k$ is
assumed to be a field of the mentioned type. An {\it affine}
monoid means a finitely generated submonoid of a free abelian
group. It is called {\it positive} if there are no non-trivial
units, and {\it simplicial} if the cone, spanned by the monoid in
the ambient Euclidean space, is simplicial (i.~e. spanned by
linearly independent vectors). An affine monoid $M$ is {\it
normal} if the implication ($c\in\NN$, $x\in\gp(M)$, $cx\in M$)\
$\Rightarrow$\ ($x\in M$) holds, $\gp(M)$ being the group of
differences of $M$. (Here we use additive notation.) $\ZZ_+$
refers to the additive monoid of non-negative integral numbers,
$\ZZ_-$, $\RR_+$ and $\RR_-$ are defined similarly and
$\NN=\{1,2,\ldots\}$. A {\it free} monoid is the one isomorphic to
$\ZZ^n_+$ for some $n\in\NN$. A cone $C\subset\RR^n$ is called
{\it unimodular} if the submonoid $C\cap \ZZ^n\subset\ZZ^n$ is
generated by a part of a basis of $\ZZ^n$. As usual, the {\it
rank} of an abelian group $H$ means $\dim_{\QQ}(\QQ\otimes H)$.
Finally, we put $r=n-1$.

\

\noindent\emph{Acknowledgement.} This paper was made possible by
the author's visit to Universit\'e Louis Pasteur (Strasbourg) and
the many discussions with Abdallah Al Amrani there. Our initial
attempt was to construct such examples among weighted projective
spaces -- a question which is motivated by the previous works
\cite{Al1,Al2} and which remains unanswered. Eventually,  I have
not been able to convince Abdallah Al Amrani to coauthor the
present paper.

Also, I am grateful to the referee for figuring out a somewhat
inaccurate use of the $\W(k)$-action on nil-$K$-theory in the
first version of the paper.

\section{Basic configuration}

Assume $M$ is a positive normal affine monoid, $\rank M=r$. (We
put $\rank M=\rank\gp(M)$.) By fixing an isomorphism in
$\gp(M)\approx\ZZ^r$ we can make the identification $\ZZ\oplus
M=\ZZ^n$. There is a free basis $\{x_1,\ldots,x_r\}$ of $\gp(M)$
such that $\RR_+M\subset\RR_+x_1+\cdots+\RR_+x_r$, both cones
being considered in $\RR^r$. In fact, passing to the dual cone
$(\RR_+M)^{\op}\subset(\RR^r)^{\op}$ we only need the existence of
a basis of the dual group $(\ZZ^r)^{\op}$ inside this dual cone
and this is easily seen (see, for instance, \cite[Lemma
2.3(d)]{G3}).

Put $e=(1,0)\in\ZZ\oplus M$ and consider the sequence of
submonoids $M=M_0$, $M_1$, $M_2$, $\ldots$ $\subset\ZZ^n$ defined
by
$M_i={\big(}\RR\oplus\RR_+M{\big)}\cap{\big(}\ZZ(x_1-ie)+\cdots+\ZZ(x_r-ie)
{\big)}$, $i\in\ZZ_+$.

\begin{lemma}\label{conf}
(a) $\ZZ\oplus M=\ZZ e +M_0=\ZZ e+M_1=\ZZ e+M_2=\cdots$, all being
direct sums,

(b) $\ZZ_+e+M_0\subset\ZZ_+e+M_1\subset\cdots$ and $\bar
M:=\bigcup_{i=0}^{\infty}{\big(}\ZZ_+e+M_i{\big)}=(\ZZ\oplus M)
\setminus\{-ie\}_{i\in\NN}$,

(c) there are isomorphisms $\alpha_i:M_i\to M_{i+1}$, $i\in\ZZ_+$,
making the diagrams
$$
\begin{CD}
\ZZ_+e+M_i@>{\subset}>>\ZZ_+e+M_{i+1}\\
@V{1+\alpha_i}VV@VV{1+\alpha_{i+1}}V\\
\ZZ_+e+M_{i+1}@>>{\subset}>\ZZ_+e+M_{i+2}
\end{CD}
$$
commutative.
\end{lemma}

This is easily proved. The isomorphisms $\alpha_i$, for instance,
are the corresponding restrictions of the automorphism
$\ZZ^n\to\ZZ^n$, $e\mapsto e$, $x_j\mapsto x_j-e$ (equivalently,
$x_j-ie\mapsto x_j-(i+1)e$), $j\in[1,r]$.

Fix arbitrarily two affine normal submonoids $N_+\subset\ZZ_+e+M$
and $N_-\subset\ZZ_-e+M_1$ so that the following conditions hold:
(i) $e\subset N_+$, $-e\subset N_-$, (ii) $N_-\cap M_1=\{0\}$, and
(iii) $\ZZ e+N_+=\ZZ e+N_-=\ZZ e+M$. This is possible: fix two
rational cones $C_+\subset\RR_+e+\RR_+M$ and
$C_-\subset\RR_-e+\RR_+M_1$ so that $C_+$ is bounded by the facets
of $\RR_+e+\RR_+M$, containing $e$, and one more hyperplane
through the origin, and similarly for the cone $C_-$ with respect
to $\RR_-e+\RR_+M_1$ under the following additional requirement:
the facet of $C_-$ that does not contain $-e$ intersects
$\RR_+M_1$ only at the origin. Then $N_+=C_+\cap\ZZ^n$ and
$N_-=C_-\cap\ZZ^n$.

A triple of the type $\C=(M,N_+,N_-)$ will be called a {\it basic
configuration}.

\section{Noncomplete case}

Let $\C$ be a basic configuration. By $X(\C)$ we denote the scheme
obtained by gluing $\Spec(k[N_+])$ and $\Spec(k[N_-])$ along their
common open subscheme $\Spec(k[\ZZ\oplus M])$. This is a toric
variety whose fan has two maximal $n$-dimensional cones -- the
dual cones
$(\RR_+N_+)^{\op},(\RR_+N_-)^{\op}\subset(\RR^n)^{\op}$. The
quasiprojecivity of $X(\C)$ is discussed in Section \ref{proj}
below.

\begin{proposition}\label{two}
Let $\C=(M,N_+,N_-)$ be a basic configuration in which the monoid
$M$ is simplicial and non-free. Then $\rank
K_0(X(\C))\geq\dim_{\QQ}k$.
\end{proposition}

We will use the fact that the $K$-groups of \cite{TT} agree with
those of Quillen for quasiprojective schemes over an affine
scheme.

\begin{proof}
We will use multiplicative notation for the monoid operation with
$X=e$. By \cite[Theorem 8.1]{TT} we have the exact sequence
$K_1(k[N_+])\oplus K_1(k[N_-])\to K_1(k[M][X,X^{-1}])\to
K_0(X(\C))$ and, hence, it is enough to show that
$$
\rank\Coker{\big(}K_1(k[N_+])\oplus K_1(k[N_-])\to
K_1(k[M][X,X^{-1}]){\big)} \geq\dim_{\QQ}k.
$$
On the other hand, $k[M][X,X^{-1}]= k[M_1][X,X^{-1}]$ by Lemma
\ref{conf}(a) and
$$
K_1(k[M_1][X,X^{-1}])=K_1(k[M_1])\oplus NK^+_1(k[M_1])\oplus
NK^-_1(k[M_1])\oplus K_0(k[M_1])
$$
by the \emph{Fundamental Theorem} \cite[Ch.7,\S7]{Ba}. Here
$K_1(k[M_1])\oplus NK^{+}_1(k[M_1])=K_1(k[M_1][X])$ and
$K_1(k[M_1])\oplus NK^{-}_1(k[M_1])=K_1(k[M_1][X^{-1}])$. Because
of the inclusions $k[N_+]\subset k[M][X]$ and $k[N_-]\subset
k+X^{-1}k[M_1][X^{-1}]$ the group
$$
\Coker{\big(}K_1(k[N_+])\oplus K_1(k[N_-])\to
K_1(k[M][\ZZ]){\big)}$$ surjects onto
$$
\Coker{\big(}K_1(k[M][X])\oplus NK_1^-(k[M_1])\to
K_1(k[M_1][X,X^{-1}]){\big)}.
$$
But the latter group contains $\Coker{\big(}K_1(k[M][X])\to
K_1(k[M_1][X]){\big)}$ as a direct summand because, firstly, the
homomorphism $k[M][X]\to k[M_1][X,X^{-1}]$ factors through
$k[M_1][X]\to k[M_1][X,X^{-1}]$ and, secondly, the groups
$K_1(k[M_1][X])$ and $NK_1^-(k[M_1])$ inside
$K_1(k[M_1][X,X^{-1}])$ are complementary direct summand. It is
therefore enough to show the inequality
\begin{align*}
\rank\Coker{\big(}K_1(k[M][X])\to K_1(k[M_1][X]){\big)}=\\
\rank\Coker{\big(}K_1(k[M][X])/k^*\to K_1(k[M_1][X])/k^*{\big)}
\geq\dim_{\QQ}k.
\end{align*}
(We view the multiplicative group $k^*$ as a natural direct
summand of the $K_1$-groups.) We can fix a grading
$k[M_1][X]=k\oplus A_1\oplus A_2\oplus\cdots$ so that all elements
of $M_1$ as well as the variable $X$ are homogeneous. This grading
restricts to a grading on $k[M][X]$.

By Weibel \cite{W}, generalizing the Bloch-Stienstra operations on
nil-$K$-theory \cite{Bl}\cite{S} to the relative situation in
arbitrary graded rings, the both groups $K_1(k[M][X])/k^*$ and
$K_1(k[M_1][X])/k^*$ are modules over the ring of big Witt vectors
$\W(k)$. The {\it ghost map} establishes a ring isomorphism
$\W(k)\approx\Pi_1^{\infty}k$. Thus there is a copy of $k$ inside
$\W(k)$. In particular, due to functoriality of the
$\W(k)$-action, the homomorphism
$$
K_1(k[M][X])/k^*\to K_1(k[M_1][X])/k^*
$$
is a homomorphism of $k$-vector spaces. Thus everything boils down
to non-sur\-jec\-ti\-vi\-ty of this homomorphism or, equivalently,
to non-surjectivity of
$$
K_1(k[M][X])\to K_1(k[M_1][X]).
$$
Assume the latter map is surjective. We have a filtered union
representation $k[\bar M]=\bigcup_i k[M_i][X]$, $\bar M$ as in
Lemma \ref{conf}(b). In particular, $K_1(k[\bar M])=
\lim_{\to}\bigl(K_1(k[M][X])\to\cdots\to
K_1(k[M_i][X])\to\cdots\bigr)$. By Lemma \ref{conf}(c) the mapping
$K_1(k[M_i][X])\to K_1(k[M_{i+1}][X])$ for every index $i$ is the
same, up to an isomorphic transformation, as the initial
homomorphism $K_1(k[M][X])\to K_1(k[M_1][X])$.  Therefore, by our
surjectivity assumtion all these mappings are surjective. In
particular, the limit map $K_1(k[M][X])\to K_1(k[\bar M])$ is also
surjective. But it is injective as well because so is the
composite map $K_1(k[M][X])\to K_1(k[\bar M])\to
K_1(k[M][X,X^{-1}])$. The diagram
$$
\begin{CD}
k[\bar M]@>{\subset}>>k[M][X,X^{-1}]\\
@VVV@VVV\\
k[X]@>>{\subset}>k[X,X^{-1}]
\end{CD}
$$
is a Cartesian square by Lemma \ref{conf}(b). Here the vertical
maps are defined by $m\mapsto0\in k$ for all non-trivial elements
$m\in M$. Using the equality $K_1(k[M][X])=K_1(k[\bar M])$ the
associated Milnor Mayer-Vietoris sequence reads as
$$
K_1(k[M][X])\to K_1(k[X])\oplus K_1(k[M][X,X^{-1}])\to
K_1(k[X,X^{-1}])
$$
(Notice, just surjectivity of $K_1(k[M][X])\to K_1(k[\bar M])$ is
enough for this sequence.) Therefore, the Fundamental Theorem
implies $NK^-_1(k[M])=0$, i.~e. $K_1(k[M])=K_1(k[M][X^{-1}])$.
Since $k[M]$ is graded ring whose 0th component is $k$ the
\emph{Swan-Weibel homotopy trick} (see below) implies
$K_1(k)=K_1(k[M])$, i.~e. $SK_1(k[M])=0$ -- a contradiction by
\cite{G2}.
\end{proof}

The mentioned homotopy trick says that for a graded ring $\Lambda=
\Lambda_0\oplus\Lambda_1\oplus\cdots$ and a functor $F$ from rings
to abelian groups the following implication holds
$F(\Lambda)=F(\Lambda[Z]) \Rightarrow F(\Lambda_0)=F(\Lambda)$,
$Z$ being an indeterminate (see \cite{An}).

\begin{corollary}\label{compl}
Let $\Zz(\Ff)$ be a quasiprojective $n$-dimensional toric variety,
defined by a fan $\Ff$ in the dual space $(\RR^n)^{\op}$, and
$\C=(M,N_+,N_-)$ be a basic configuration in which the monoid $M$
is simplicial and nonfree. Assume the dual cones
$(\RR_+N_+)^{\op},(\RR_+N_-)^{\op}\subset(\RR^n)^{\op}$ are among
the maximal cones of $\Ff$ and all other maximal cones of $\Ff$
are unimodular. Then $\rank K_0(\Zz)\geq\dim_{\QQ}k$.
\end{corollary}

\begin{proof}
We have the open cover $\Zz=X(\C)\cup\Aa_1\cup\cdots\cup\Aa_s$,
where the $\Aa_j$ are the \emph{smooth} affine toric varieties
(i.~e. of type $\AA_k^a\times\TT_k^b$, $a+b=n$) that correspond to
the mentioned unimodular maximal cones. By \cite[Theorem 8.1]{TT}
for each $j\in[1,s]$ we get the exact sequence
\begin{multline*}
K_0{\big(}X(\C)\cup\Aa_1\cup\cdots\cup\Aa_j{\big)}\to
K_0{\big(}X(\C)\cup\Aa_1\cup\cdots\cup\Aa_{j-1}{\big)}\oplus
K_0(\Aa_j)\to K_0(U_j),
\end{multline*}
where $U_j={\big(}X(\C)\cup\Aa_1\cup\cdots\cup\Aa_{j-1}{\big)}
\cap\Aa_j$. Since the $U_j$ are open subschemes of $\Aa_j$ and
$K_0(\Aa_j)=\ZZ$ we have $K_0(U_j)=\ZZ$, $j\in[1,s]$. Therefore,
the induction from $j=1$ to $j=s$ shows that $\rank
K_0(\Zz)\geq\dim_{\QQ}k$. (Here one needs $\dim_{\QQ}k=\infty$.)
\end{proof}

\section{Projective case}\label{proj}

Given a basic configuration $\C$ in which the monoid $M$ is
simplicial. Then $X(\C)$ can be embedded into a projective
simplicial toric variety as an equivariant open subscheme. A quick
way to do so is given by the following construction. Consider the
intersection $\Delta(\C)=(e+\RR_+N_-)\cap\RR_+N_+$. This is a
rational $n$-simplex having a pair of corners, spanning
respectively the cones $\RR_+N_+$ and $e+\RR_+N_-$. Let
$\overline{X(\C)}$ be the following projective toric variety. We
have the $(n+1)$-dimensional cone
$C(\Delta(\C))=\bigcup_{x\in\Delta(\C)}
\RR_+(x,1)\subset\RR^{n+1}$ and the associated monoid ring
$k[C(\Delta(\C))\cap\ZZ^{n+1}]$, which carries the graded
structure $k[C(\Delta(\C))\cap\ZZ^{n+1}]=k\oplus C_1\oplus
C_2\oplus\cdots$ given by the last coordinate in the exponent
vectors of monomials. Then
$\overline{X(\C)}=\Proj{\big(}k[C(\Delta(\C))\cap\ZZ^{n+1}]{\big)}$.
Its standard affine charts are $\Spec(k[C_v\cap\ZZ^n])$,
$v\in\vert(\Delta(\C))$, where the $C_v\subset\RR^n$ are the
corner cones spanned by $\Delta(\C)$ at its vertices $v$ and then
shifted by $-v$. The maximal cones in the fan of
$\overline{\X(\C)}$ are just the dual cones
$(-v+C_v)^{\op}\subset(\RR^n)^{\op}$, $v\in\vert(\Delta(\C))$
\cite[Section 1.5]{F}. This fan contains $(\RR_+N_+)^{\op}$ and
$(\RR_+N_-)^{\op}$ as adjacent maximal cones whose common facet is
exactly $(\RR e+\RR_+M)^{\op}$.

The difficulty in applying Corollary \ref{compl} to
$\overline{X(\C)}$ is that the cones $-v+C_v$,
$v\in\vert(\Delta(\C))$, different from $\RR_+N_+$ and $\RR_+N_-$,
are in general not unimodular. We overcome this difficulty by
resolving the corresponding toric singularities without affecting
$\Spec(k[N_+])$ and $\Spec(k[N_-])$. As we will see, sometimes
this is possible.

Call a basic configuration $\C=(M,N_+,N_-)$ {\it admissible} if
$M$ is simplicial and all facets of $(\RR_+N_+)^{\op}$ and
$(\RR_+N_-)^{\op}$, except maybe their common facet, are
unimodular.

\begin{lemma}\label{adm}
Assume $\C$ is an admissible basic configuration. Then there
exists an equivariant open embedding of $X(\C)$ into a projective
simplicial toric variety $\Zz$ whose affine charts that do not
come from $X(\C)$ are all smooth.
\end{lemma}

\begin{proof}
We apply the standard equivariant resolution to $\overline{X(\C)}$
\cite[Section 2.6]{F} except we do not touch the cones
$(\RR_+N_+)^{\op},(\RR_+N_-)^{\op}$. This is possible because the
admissibility condition on $\C$ guarantees that the
cone-subdividing elements $z\in(\ZZ^n)^{\op}$ we produce in the
resolution process do not belong to the boundary
$\partial{\big(}(\RR_+N_+)^{\op}\cup(\RR_+N_-)^{\op}{\big)}$, that
is
$z\in(\RR^n)^{\op}\setminus{\big(}(\RR_+N_+)^{\op}\cup(\RR_+N_-)^{\op}{\big)}$.
The projectivity condition is not lost in this process because we
essentially have barycentric (rather, \emph{stellar}) subdivisions
-- they are projective -- and a composition of projective
subdivisions is again projective \cite[Ch.3, \S1]{KKMS}.
\end{proof}

Now we are ready to prove

\begin{theorem}\label{main}
For each $n\geq3$ there are projective simplicial toric varieties
$\Zz$ for which $\rank K_0(\Zz)\geq\dim_{\QQ}k$.
\end{theorem}

\begin{proof}
In view of Corollary \ref{compl} and Lemma \ref{adm} we only need
to show the existence of admissible configurations
$\C=(M,N_+,N_-)$ in which $M$ is not free.

First we observe that there are simplicial rational $n$-cones
whose all facets except one are unimodular and the distinguished
facet is not unimodular (see Example \ref{ex} below). Fix such a
cone $C$ in the dual space $(\RR^n)^{\op}$ (with respect to the
dual lattice $(\ZZ^n)^{\op}$). Let $F\subset C$ be the
nonunimodular facet. The dual cone $C^{\op}$ is a rational
simplicial cone whose one edge, say $l$, corresponds to $F$ (under
the duality). Put $N_+=C^{\op}\cap\ZZ^n$ and denote by $e$ the
generator of $\ZZ^n\cap l\approx\ZZ_+$. By \cite[Theorem 1.8]{G1}
$\ZZ_-e+N_+=\ZZ e + M'$ for some rank $r$ simplicial normal monoid
$M'$. We have $(\RR e+\RR_+M')^{\op}=F$. As in Lemma \ref{conf} we
can find a free basis $\{x_1,\ldots,x_r\}$ of $\gp(M')$ such that
$\RR_+M'\subset\RR_+x_1+\cdots+\RR_+x_r$. If a natural number $c$
is big enough then the monoid
$M=\ZZ^n\cap{\big(}\RR_+(x_1-ce)+\cdots+\RR_+(x_r-ce){\big)}$
satisfies the conditions $\RR_+N_+\subset\RR_+e+\RR_+M$
(equivalently, $N_+\subset\ZZ_+e+M$) and $\RR_+N_+\cap\RR_+M=0$.
We fix such $c$ and define the monoid $M_1$ to be
$\ZZ^n\cap{\big(}\RR_+(x_1-(c+1)e)+\cdots+\RR_+(x_r-(c+1)e){\big)}$.
To construct $N_-$ we make the identification of the monoids
$\ZZ_+e+M$ and $\ZZ_-e+M_1$ along the isomorphism induced by
$e\mapsto -e$, $x_j-ce\mapsto x_j-(c+1)e$ ($j\in[1,r]$). Then
$N_-$ is by definition the corresponding copy of $N_+$. It follows
from the construction that $(M,N_+,N_-)$ is an admissible basic
configuration. (Say, the non-freeness of $M$ follows from the
non-unimodularity of $F$, i.~e. non-smoothness of $\Spec(k[\ZZ
e+M])$.)
\end{proof}

\begin{example}\label{ex}
For each natural number $r\geq2$ consider the simplicial normal
nonfree monoid
$M_r=\ZZ_+(1,\ldots,1)+\sum_{j=1}^r\ZZ_+(re_j)\subset\ZZ_+^r$,
where $e_j$ is the $j$th standard basic vector in $\RR^r$
($j\in[1,r]$). The facets of $\RR_+^r=\RR_+M_r$ are unimodular
with respect to the sublattice $\gp(M_r)\subset\ZZ^r$. Therefore,
all facets of the cone $\RR_+(\ZZ_+\oplus M_r)=\RR^n_+$, except
exactly one, are unimodular with respect to the lattice
$\gp(\ZZ_+\oplus M_r)=\ZZ\oplus\gp(M_r)\subset\RR^n$.
\end{example}


\begin{thebibliography}{******}
\bibitem[Al1]{Al1}
A. Al-Amrani, {\em Groupe de Chow et $K$-th\'eorie coh\'erente des
espaces projectifs tordus}, $K$-Theory {\bf 2} (1989), 579--602.
\bibitem[Al2]{Al2}
A. Al-Amrani, {\em Complex $K$-theory of weighted projective
spaces}, J. Pure Appl. Algebra {\bf 93} (1994), 113--127.
\bibitem[An]{An}
D. F. Anderson, {\em Projective modules over subrings of
$k[X,\,Y]$ generated by monomials},  Pacific J. Math. {\bf 79}
(1978), 5--17.
\bibitem[Bl]{Bl}
S. Bloch, {\em Some formulas pertaining to the K-theory of
commmutative groupschemes}, J. Algebra {\bf 53} (1978), 304--326.
\bibitem[Ba]{Ba}
H. Bass, {\em Algebraic $K$-theory}, W. A. Benjamin, Inc., 1968.
\bibitem[BFM]{BFM}
P. Baum, W. Fulton, R. MacPherson, {\em Riemann-Roch for singular
varieties}, Inst. Hautes \'Etudes Sci. Publ. Math. No. {\bf 45}
(1975), 101--145.
\bibitem[BV]{BV}
M. Brion, M. Vergne, {\em An equivariant Riemann-Roch theorem for
complete, simplicial toric varieties}, J. Reine Angew. Math. {\bf
482} (1997), 67--92.
\bibitem[C]{C} D. A. Cox, {\em Update on toric geometry},
S\'eminaires et Congr\`es {\bf 6} 2002, 1--42.
\bibitem[F]{F}
W. Fulton, {\em  Introduction to toric varieties}, Princeton
University Press, 1993.
\bibitem[FMSS]{FMSS}
W. Fulton, R. MacPherson, F. Sottile, B. Sturmfels, {\em
Intersection theory on spherical varieties}. J. Algebraic Geom.
{\bf 4} (1995), 181--193.
\bibitem[G1]{G1}
J. Gubeladze, {\em Anderson's conjecture and the maximal class of
monoids over which projective modules are free}, Math. USSR
Sbornik {\bf 135} (1988), 169--181 (Russian).
\bibitem[G2]{G2}
J. Gubeladze, {\em Nontriviality of $SK_1(R[M])$,} J. Pure Appl.
Algebra {\bf 104} (1995), 169--190.
\bibitem[G3]{G3}
J. Gubeladze, {\em Higher $K$-theory of toric varieties},
$K$-Theory, to appear.
\newline(http://xxx.lanl.gov/abs/math.KT/0104166)
\bibitem[KKMS]{KKMS}
G. Kempf, F. Knudsen, D. Mumford, and B. Saint-Donat, {\em
Toroidal embeddings~I}, Lecture Notes in Math. {\bf 339},
Springer, 1973.
\bibitem[S]{S}
J. Stienstra, {\em Operations in the higher K-theory of
endomorphisms. (Current trends in algebraic topology, Semin.
London/Ont. 1981)}, CMS Conf. Proc. 2, Amer. Math. Soc.,
Providence, R.I., 1982, 59--115.
\bibitem[TT]{TT}
R. W. Thomason and T. F. Trobaugh, {\em Higher algebraic
$K$-theory of Schemes and of derived categories}, The Grothendieck
Festschrift III, Prog. Math. {\bf 88}, Birkh\"auser, 1990,
247--435.
\bibitem[W]{W}
C. Weibel, {\em Module structures on the $K$-theory of graded
rings}, J. Algebra {\bf 105} (1987), 465--483
\end{thebibliography}
\end{document}